\documentclass[12pt,article]{amsart}

\usepackage[margin=1in]{geometry}
\usepackage{amsmath}
\usepackage{cancel}
\usepackage{amssymb}
\usepackage{epsfig}
\usepackage[normalem]{ulem}
\usepackage[usenames,dvipsnames]{color}
\usepackage[textwidth=20mm]{todonotes}
\usepackage{hyperref}
\usepackage{eqnarray}
\usepackage{mathtools}
\usepackage{mathrsfs}

\usepackage{algpseudocode}
\usepackage{algorithm}
\usepackage{tcolorbox}
\usepackage{multicol}
\usepackage{epigraph}
\usepackage{enumerate}

\usepackage{graphicx}
\graphicspath{ {./images/} }

\usepackage{soul}
\usepackage{tikz}

\usetikzlibrary{backgrounds}
\usetikzlibrary{intersections}

\usepackage{float}
\usepackage{xcolor}
\colorlet{shadecolor}{gray!15}

\newtheorem{theorem}{Theorem}[section]

\newtheorem{lemma}[theorem]{Lemma}
\newtheorem{cor}[theorem]{Corollary}

\newtheorem{conj}{Conjecture}

\newcommand{\fancyP}{\mathcal{P}}

\newcommand{\E}{\mathbb E}

\newcommand{\Pro}{\mathbb{P}}

\usepackage{framed}

\title{Almost all 9-regular graphs have a modulo-5 orientation}

\author{Michelle Delcourt}
\address{Department of Mathematics, Toronto Metropolitan University (formerly named Ryerson University), Toronto, Ontario M5B 2K3, Canada}
\email{\tt mdelcourt@torontomu.ca}
\thanks{Research supported by NSERC under Discovery Grant No.\ 2019-04269.}

\author{Reaz Huq}
\address{Department of Mathematics, Toronto Metropolitan University (formerly named Ryerson University), Toronto, Ontario M5B 2K3, Canada}
\email{\tt reaz.huq@torontomu.ca}

\author{Pawe\l{} Pra\l at}
\address{Department of Mathematics, Toronto Metropolitan University (formerly named Ryerson University), Toronto, Ontario M5B 2K3, Canada}
\email{\tt pralat@torontomu.ca}
\thanks{Research supported by NSERC under Discovery Grant No.\ 2022-03804. Part of this work was done while the author was visiting the Simons Institute for the Theory of Computing.}

\date{\today}

\begin{document}

\maketitle

\begin{abstract}
In 1972 Tutte famously conjectured that every 4-edge-connected graph has a nowhere-zero 3-flow; this is known to be equivalent to every 5-regular, 4-edge-connected graph having an edge orientation in which every in-degree is either 1 or 4.  In 1988 Jaeger conjectured a generalization of Tutte's nowhere-zero 3-flow conjecture, by stating that every $(4p+1)$-regular, $4p$-edge-connected graph has an edge orientation in which every in-degree is either $p$ or $3p+1$.  Inspired by the work of Pra\l{}at and Wormald investigating $p=1$, we address $p=2$ to show that the conjecture holds asymptotically almost surely for random 9-regular graphs. It follows that the conjecture holds for almost all 9-regular, 8-edge-connected  graphs.  These results make use of the technical small subgraph conditioning method.
\end{abstract}

\section{Introduction}

A \emph{$k$-flow} of an undirected graph $G=(V,E)$ is an orientation of its edge-set together with a function $f: E \rightarrow \{0,\pm 1,\pm 2, \ldots, \pm(k-1)\}$ such that the following is satisfied for each vertex $v \in V$:
$$
\sum\limits_{e \in D^+(v)}f(e) - \sum\limits_{e \in D^-(v)}f(e) = 0.
$$
To orient a graph $G=(V,E)$, convert each edge $e=\{x,y\}$ into either the ordered pair $e=(x,y)$ or $e=(y,x)$, to represent the edge being oriented from $x$ to $y$ or from $y$ to $x$, respectively. We denote the set of all edges oriented towards $v$ as $D^{+}(v)$ and the set of all edges oriented away from $v$ as $D^{-}(v)$; a $k$-flow is said to be a \emph{nowhere-zero $k$-flow} if $f(e)\neq 0$ for all $e \in E$. We define the \emph{in-degree} of $v$ as $d^{+}(v):=|D^{+}(v)|$ and the \emph{out-degree} of $v$ as $d^{-}(v):=|D^{-}(v)|$. In 1972 Tutte famously conjectured the following (for instance, see Bondy \cite[Open Problem 48]{bondy} or Jensen \cite[Section 13.3]{jensen}):
\begin{conj}(Tutte's nowhere-zero 3-flow conjecture)
Every 4-edge-connected graph admits a nowhere-zero 3-flow.
\end{conj}
Though Tutte's nowhere-zero 3-flow conjecture has attracted considerable attention, it has yet to be proven. For a significant period of time it was not known whether or not there exists a fixed constant $k$ such that every $k$-edge connected graph has a nowhere-zero 3-flow (this is known as the weak 3-flow conjecture of Jaeger). It was proved for $k \geq c \log_2 n$ for $n$-vertex graphs by Alon, Linial, and Meshulam in 1991 \cite{nogalm}, as well as by Lai and Zhang in 1992 \cite{laiz}. The weak $3$-flow conjecture was first proved in 2012 by Thomassen \cite{thomas}, who showed that every $8$-edge-connected graph admits a nowhere-zero 3-flow. This was improved to $k=6$ in 2011 by Lov\'asz, Thomassen, Wu, and Zhang \cite{ltwz}. It is a well-known fact that a graph admits a nowhere-zero 3-flow if and only if it has a nowhere-zero 3-flow over $\mathbb{Z}_3$, for instance see Seymour \cite{seymour} for more discussion.  The following is a reformulation of Tutte's nowhere-zero 3-flow conjecture (\cite[Section 13.3]{jensen}, \cite[Open Problem 48]{bondy}):
 \begin{conj}(Tutte's nowhere-zero 3-flow conjecture reformulated)
Every 5-regular, 4-edge-connected graph has an edge orientation in which
every in-degree is either 1 or 4.
\end{conj}

Let $k \geq 3$ be an odd integer; an orientation is a \emph{modulo k-orientation} if for every vertex $v \in V$, 
$d^{+}(v) \equiv d^{-}(v) \pmod{k}$. In 1988 Jaeger \cite{jaeger} conjectured the following for all $p$ of which Tutte's nowhere-zero 3-flow conjecture is a sub-case when $p=1$:

\begin{conj}{(Jaeger's conjecture)}
For any integer $p \geq 1$, the edges of every $(4p+1)$-regular, $4p$-edge-connected graph has a modulo $(2p + 1)$-orientation, i.e.\ can be oriented so that every in-degree is either $p$ or $3p+1$. 
\end{conj}

\noindent Jaeger's conjecture was originally formulated for $4p$-edge-connected graphs but it can be reduced to the case of $(4p + 1)$-regular graphs (see for instance \cite{pralatflow}, \cite{delc}, and \cite{alon}, which use the same formulation of the conjecture, and the following from Jaeger: \cite[Theorem 5.9]{jaeger}, \cite[Remark 2]{jaeger1984}). In 2019, Pra\l at and Wormald \cite{pralatflow} used the small subgraph conditioning method to demonstrate that Jaeger's conjecture is asymptotically almost surely true for $p=1$. In 2011, Alon and Pra\l at \cite{alon} used the expander mixing lemma to demonstrate that Jaeger's conjecture is asymptotically almost surely true for $p$ which are somewhat large (though still finite). In 2018, Han, Li, Wu, and Zhang \cite{HLWZ18} demonstrated the existence of $4p$-edge-connected graphs which do not admit a $(2p+1)$-orientation for $p \geq 3$.  A natural question is whether or not the conjecture holds asymptotically almost surely for the values of $p$ in between the value of $p=1$ identified by Pra\l at and Wormald \cite{pralatflow} and the large -- though not explicitly determined -- values of $p$ identified by Alon and Pra\l at \cite{alon}. 

We answer the question in the affirmative in the case where $p=2$. Our results are asymptotic, that is, we are holding that the probability that a random graph on $n$ vertices has a particular property tends to 1 as $n$ tends to infinity. By utilizing the pairing model of Bollob\'as (see \cite{bolo}) -- which a.a.s.\ produces 9-regular, 8-edge connected graphs -- we investigate the case where $p=2$ and study $9$-regular, 8-edge connected graphs to prove the following theorem: 
\begin{theorem}\label{thm:jaegp}
A random 9-regular graph $G_n$ on $n$ vertices a.a.s. has a modulo-5 orientation, i.e. has an orientation in which every in-degree is either 2 or 7. 
\end{theorem} 

Our result is of relevance to conjectures besides that of Jaeger's. Modular 5-orientations in highly connected graphs are particularly interesting object.  As Han, Li, Wu, and Zhang noted in~\cite{HLWZ18}, if $8$-edge-connected graphs in general admit modular 5-orientations, it is implied that Tutte's famous 5-Flow Conjecture holds.
 
\section{Overview of Our Proof of the Main Result}

To prove Theorem~\ref{thm:jaegp} we utilize the famous pairing model introduced by  Bollob\'as (for more details, see \cite{bolo}).  In the pairing model we work with a random graph $\mathcal{P}_{n,d}$, in which $dn$ points are arranged in $n$ groups of $d$ points each. Each of the $n$ groups is referred to as a \emph{vertex} while the $d$ points composing them continue to be referred to as \emph{points}. A perfect matching is then made amongst the points with each group of $d$ points forming a vertex of degree $d$. The result is a $d$-regular multi-graph. The convenience of the pairing model lies in the fact that, when $d = O(1)$, if a given statement is asymptotically almost surely true in $\mathcal{P}_{n,d}$ then it is asymptotically almost surely true in the random $d$-regular graph on $n$ vertices (see Corollary 2.3 in \cite{wormald}).

We are interested in orientations of the edge-set of 9-regular graphs in which every vertex has an in-degree of either 2 or 7 (which implies an out-degree of either 7 or 2, respectively); we call such orientations \emph{valid}. An orientation of an element of $\fancyP_{n,9}$ is referred to as ``valid" when each of the $n$ groups of $d$ points has either exactly two of its points serving as the terminal ends of two edges after orienting them, or exactly two of its points serving as the initial ends of two edges after orienting them.  Let $Y=Y(n)$ be the number of valid orientations of a random element of $\fancyP_{n,9}$. 

In Section~\ref{sec:exp}, we approximate $\E [Y]$ using Stirling's approximation:

\begin{lemma}\label{lemma:exp}
$$\E[Y]  \sim 3 \left(\frac{81}{8} \right)^{n/2} > 0.$$
\end{lemma}

We will show that there are plenty of valid orientations per pairing, on average. In order to show that pairings a.a.s.\ have at least one valid orientation (that is, $\Pro(Y > 0) \sim 1$ or, alternatively $\Pro(Y = 0) \sim 0$), a~standard strategy is to estimate $\E [Y(Y-1)]$ and show that it is asymptotic to $(\E [Y])^2$.  To this end, in Section~\ref{sec:second} using optimization, Taylor series expansions
and multivariate calculus we show the following approximation:

\begin{lemma}\label{lemma:second}
$$\E[Y(Y-1)] \sim \left( \frac{81}{8} \right)^n \frac{81}{7}.$$
\end{lemma}

One could hope that at this point the desired result would then immediately follow from Chebyshev's inequality:
$$
\Pro \left( Y = 0 \right) \le \frac { \E [Y(Y-1)] }{ (\E [Y])^2 } - 1. 
$$
Unfortunately, this proof technique fails here as there is a constant factor discrepancy in the asymptotic ratio giving us an upper bound of $2/7+o(1)$ for the failure probability instead of $o(1)$. 

\begin{cor}\label{cor:second}
$$
\frac{\E[Y(Y-1)]}{\E[Y]^2} \sim \frac{9}{7}>0.
$$
\end{cor}

Instead to prove our main result Theorem~\ref{thm:jaegp}, we will utilize the powerful but highly technical small subgraph conditioning method introduced by Robinson and Wormald in 1994~\cite{RW94}.  We refer the curious reader to a survey paper by Wormald~\cite{wormald} for a more detailed discussion about the small subgraph conditioning method.  This technique was used by Pra\l at and Wormald to demonstrate the truth of Jaeger's conjecture for $p = 1$ \cite{pralatflow} and it has been used to deal with similar problems, as in \cite{delc}, where Delcourt and Postle used it to demonstrate that a 4-regular graph $G$ a.a.s.\ admits an $S_3$-decomposition, that is, $G$ can be shown to consist of copies of three vertices which all make an edge with another vertex. This line of inquiry in utilizing the small subgraph conditioning method to study star-decompositions was extended by Delcourt, Greenhill, Isaev, Lidick{\'y}, and Postle in \cite{delgreen}, where it is shown that the uniform random $d-$regular graph admits a $k-$star decomposition for all $d,k$ such that $d/2 < k \leq d/2 + \max\{1,\log d /6 \}$. We note that $[x]_k :=x(x-1)\cdots(x-k+1)$ denotes the falling factorial:

\begin{theorem}[Small Subgraph Conditioning Method]\label{thm:sscm}
Let $\lambda_j>0$ and $\delta_j>-1$ be real numbers for all $j \geq 0$.  Suppose
for each $n$ there are non-negative random variables $X_j = X_j(n)$, $j \geq 1$, and $Y = Y (n)$ (defined on the same probability space) such that $X_j$ is integer valued and $\E[Y ] > 0$ (for $n$ sufficiently large).
Furthermore, suppose that
\begin{enumerate}[(a)]\label{sscm}
\item For each $k \geq 1$, $X_i$ ($i \in \{1,2,\ldots,k\}$) are asymptotically independent Poisson random variables with $\E[X_i] \sim \lambda_i$;
\item For every finite sequence $j_1, j_2, \ldots, j_k$ of non-negative integers,
$$
\frac{\E[Y[X_1]_{j_1} \ldots [X_k]_{j_k}]}{\E [Y]} \sim \prod\limits_{i=1}^k \mu_i^{j_i} = \prod\limits_{i=1}^k \left( \lambda_i (1+ \delta_i) \right)^{j_i};
$$
\item $\sum_i \lambda_i \delta_i^2 < \infty $;
\item $$
\frac{\E [Y^2]}{(\E [Y])^2} \leq \exp \left( \sum\limits_i \lambda_i \delta_i^2 \right) + o(1).
$$
\end{enumerate}
If the above four properties (a)-(d) hold, then
$$
\mathbb{P}(Y >0) = \exp \left( - \sum\limits_{\delta_i = - 1 } \lambda_i \right) + o(1).
$$
\end{theorem}

As with most applications of this method in the literature, we let $X_j$ denote the number of cycles of length $j$ in the multi-graph corresponding to a random element of $\mathcal{P}_{n,9}$.  Informally, we are in a situation where the distribution of $Y$ is affected by small but not too common (i.e., the expected number is bounded) subgraphs of the random 9-regular graph, namely short cycles, and by conditioning on the small subgraph counts, we are able to control the variance of $Y$ and show that $Y > 0$ asymptotically almost surely.  In our setting, condition (a) of the Small Subgraph Conditioning Method (Theorem~\ref{thm:sscm}) follows from a well known result of Bollob\'as~\cite{bolo}, namely that for $j\geq1$ $X_1,X_2, \ldots, X_j$ are asymptotically independent Poisson distributions with mean
$$
\E [X_j] \longrightarrow \lambda_j := \frac{8^j}{2j}.
$$

For condition (b) of Theorem~\ref{thm:sscm}, in Section~\ref{sec:joint}, we compute  $\E [Y X_j]$ and obtain that
\begin{lemma}\label{lemma:cycles}
For all $j \geq 1$
$$\frac{\E [Y X_j ]}{\E Y} = \frac{8^j}{2j} \left( 1 + \left( -\frac{2}{9} \right)^j \right) = \lambda_j \left( 1 + \left( -\frac{2}{9} \right)^j \right).$$
\end{lemma}
\noindent
Therefore, for all $j \geq 1$
$$\delta_j := -\left(\frac{2}{9}\right)^j > -1.$$

For any fixed subgraph $H$ which has more edges than vertices, a multigraph corresponding to a random element of $\mathcal{P}_{n,9}$ a.a.s.\ contains no subgraph isomorphic to $H$. As in other applications of the Small Subgraph Conditioning Method, we would not expect two cycles to sharing edges (or for that matter vertices), and therefore, we concentrate on disjoint cycles and think of them roughly as being independent. These observations combined with Lemma~\ref{lemma:cycles} imply the following more general form of Lemma~\ref{lemma:cycles}, computing joint
factorial moments:
\begin{cor}
For all $j \geq 1$, if $(\ell_1, \ell_2, \ldots, \ell_j)$ is a sequence of non-negative integers, the following holds:
$$
\frac{\E(Y [X_1]_{\ell_1} \cdots [X_j]_{\ell_j})}{\E Y} \sim \prod\limits_{i=1}^j \mu_i^{\ell_i} = \prod\limits_{i=1}^j \left(  \lambda_i(1+\delta_i ) \right)^{\ell_i}
= \prod\limits_{i=1}^j \left(  \frac{8^i}{2i}\left(1 - \left(\frac{2}{9}\right)^i \right) \right)^{\ell_i}$$
\end{cor}

We now need to verify conditions (c) and (d) of Theorem \ref{thm:sscm}:
\begin{lemma}\label{lemma:condc}
$$\sum_i \lambda_i \delta_i^2 < \infty \text{ and }\exp\left(\sum_i \lambda_i \delta_i^2\right) = \frac{9}{7} \sim \frac{\E (Y^2)}{\E[Y]^2}.$$
\end{lemma}

\begin{proof}
Recall that $\lambda_k = \frac{8^k}{2k}$ and $\delta_k = -\left(\frac{2}{9}\right)^k$.  Using the fact that $-\log(1-x) = \sum_{k \geq 1} x^k/k$, 
$$
\exp\left( \sum\limits_{k \geq 1} \lambda_k \delta_k^2 \right) = \exp \left( \sum\limits_{k \geq 1} \frac{8^k}{2k} \left( \frac{4}{81} \right)^k \right) = \exp\left(\frac{1}{2} - \log\left( 1 - \frac{2^5}{3^4} \right) \right) = \frac{9}{7}. 
$$
Since $\exp\left( \sum\limits_{k \geq 1} \lambda_k \delta_k^2 \right) \sim \frac{\E (Y^2)}{\E[Y]^2} \sim 9/7.$ 
\end{proof}

Putting this all together, assuming the proofs of Lemma~\ref{lemma:exp} (Section~\ref{sec:exp}), Lemma~\ref{lemma:second} (Section~\ref{sec:second}), and Lemma~\ref{lemma:cycles} (Section~\ref{sec:joint}) we are now ready to prove our main result via the Small Subgraph Conditioning Method as follows:
\begin{proof}[Proof of Main Result (Theorem~\ref{thm:jaegp}).]
Let 9 divide $n$.  Let $Y=Y(n)$ be the number of valid orientations (orientations of the edge-set in which every vertex has an in-degree of either 2 or 7) of a random element of $\fancyP_{n,9}$ and let $X_j$ denote the number of cycles of length $j$ in a random element of $\fancyP_{n,9}$.  We apply the Small Subgraph Conditioning Method (Theorem~\ref{thm:sscm}) to $Y$ and the $X_j$.  Similarly to most applications in the literature, condition (a) holds due to a well-known result of Bollob\'as~\cite{bolo}, here with $\lambda_j := \frac{8^j}{2j}$.  Condition (b) holds by Lemma~\ref{lemma:cycles} and $\delta_j:=-(\frac{2}{9})^j$.  Condition (c) holds by Lemma~\ref{lemma:condc}.  Condition (d) follows from Lemma~\ref{lemma:condc} and Lemma~\ref{lemma:second}.  Hence, all four conditions of the small subgraph conditioning method are satisfied, and we may conclude that $\mathbb{P}(Y > 0) \sim 1$, as desired.
\end{proof}

\section{Expected Number of Valid Configurations}\label{sec:exp}
Recall that the random variable $Y = Y(n)$ was introduced to measure the number of valid orientations of an element of $\mathcal{P}_{n,9}$. A vertex in a valid orientation is called an \textit{in-vertex} if it has an in-degree of 2 and is called an \textit{out-vertex} otherwise. Each point contained in an edge oriented towards an in-vertex or away from an out-vertex is called \textit{special}. A point is referred to as an \textit{in-point} if the edge containing it is pointing towards it and is referred to as an \textit{out-point} otherwise.  We are now prepared to prove Lemma~\ref{lemma:exp} as follows. 

\begin{proof}[Proof of Lemma~\ref{lemma:exp}]
The expected number of valid orientations is computed as follows:
\begin{equation}\label{exp}
\E[Y] = \frac{{n \choose n/2}{9 \choose 2} ^n (9n/2)! }{M(9n)},
\end{equation}
where 
$$
M(9n) = \frac{(9n)!}{(9n/2)! 2^{9n/2}}
$$
is the number of perfect matchings on $9n$ points. Indeed, since the number of in-vertices must match the number of out-vertices, there are ${n \choose n/2}$ ways to choose in-vertices. There are ${9 \choose 2} ^n$ different ways to pick special points. Finally, there are $(9n/2)!$ different ways to pair in-points with out-points. After expanding the expression and using Stirling's formula ($x! \sim \sqrt{2\pi x} (x/e)^x$) we get that
$$
\E[Y] = \frac{ n! \, (9 \cdot 4)^n \, (9n/2)!^2 \, 2^{9n/2} }{ (n/2)!^2 \, (9n)! } \sim 3 \left(\frac{81}{8} \right)^{n/2}.
$$
\end{proof}

\section{The Second Moment Method}\label{sec:second}

To verify Lemma~\ref{lemma:second}, we now turn our attention to estimating $\E[Y(Y-1)]$. Given two orientations of the same graph, two vertices $v$ and $w$ can be related in one of three ways:
\begin{enumerate}[(i)]
\item $v$ and $w$ are in-vertices in both orientations;
\item $v$ and $w$ are out-vertices in both orientations;
\item One of $v$ or $w$ is an in-vertex in one orientation and the other is an out-vertex in the other orientation.
\end{enumerate}

Observe that if exactly $k$ vertices are in-vertices in both orientations, since there are in total $n/2$ in-vertices in the first orientation, $n/2-k$ vertices are in-vertices in the first orientation and out-vertices in the second orientation. Since there are $n/2 - k$ vertices which are out-vertices in the second orientation but not in the first and there are $n/2$ total out-vertices in the second orientation, $k$ vertices are out-vertices in both orientations. The remaining $n/2 -k$ vertices are out-vertices in the second orientation and in-vertices in the second.

We use the term $k_{211}$ to denote the number of vertices which are in-vertices in both orientations and which have two special points in common. We use the term $k_{111}$ to denote the number of vertices which are in-vertices in both orientations but which have only one special point in common. Note that $k_{211} \leq k$ and $k_{111} \leq k$ since $k_{211} + k_{111} \leq k$. It follows that there are $k-k_{211} - k_{111}$ vertices which are in-vertices in both orientations but which have no special points in common. 

Similarly, we use the term $k_{200}$ to denote the number of vertices which are out-vertices in both orientations and have two special points in common. $k_{100}$ is defined as the number of vertices which are out-vertices in both orientations but which only have one special point in common. Note that $k_{200} \leq k$ and $k_{100} \leq k$ since $k_{200} + k_{100} \leq k$. It follows that there are $k-k_{200} - k_{100}$ vertices which are out-vertices in both orientations but which have no special points in common. 

The term $k_{210}$ is used to denote the number of vertices which are in-vertices in the first orientation but out-vertices in the second and which have both special points in common. $k_{110}$ denotes the number of points which are in-vertices in the first orientation and out-vertices in the second but which have only one special point in common. Note that $k_{210} \leq n/2 - k$ and $k_{110} \leq n/2 - k$ since $k_{210} + k_{110} \leq n/2 - k$. It follows that there are $n/2 - k - k_{210} - k_{110} $ vertices which are in-vertices in the first orientation and out-vertices in the second orientation but which have no special points in common.

The term $k_{201}$ is used to denote the number of vertices which are out-vertices in the first orientation but in-vertices in the second and which have both special points in common. $k_{101}$ denotes the number of vertices which are out-vertices in the first orientation and in-vertices in the second but which have only one special point in common. Note that $k_{201} \leq n/2 - k$ and $k_{101} \leq n/2 - k$ since $k_{201} + k_{101} \leq n/2 - k$. It follows that there are $n/2 - k - k_{201} - k_{101} $ vertices which are out-vertices in the first orientation and in-vertices in the second orientation but which have no special points in common.

\subsection{Multivariate Calculus}\label{subsec:multi}
To approximate the second moment we will need to define a region that contains every possible configuration of two orientations. Let
\begin{align*}
I := I(n) = \{ & (k,k_{211},k_{111},k_{210},k_{110},k_{201},k_{101},k_{200},k_{100}) \in \mathbb{N}_{0}^{9} \bigm|  k \leq n/2, k_{211} + k_{111} \leq k,\\
&k_{200} + k_{100} \leq k, k_{210} + k_{110} \leq n/2 - k, k_{201} + k_{101} \leq n/2 - k \},
\end{align*}
where $\mathbb{N}_0 = \mathbb{N} \cup \{0\}$. Fixing $\textbf{k} = (k,k_{211},k_{111},k_{210},k_{110},k_{201},k_{101},k_{200},k_{100})$, we can calculate the number of configurations corresponding to this vector. Letting $k_{011} = k-k_{211}-k_{111}$, $k_{000} = k-k_{200}-k_{100}$, $k_{010} = n/2-k-k_{210}-k_{110}$, and $k_{001} = n/2-k-k_{201}-k_{101}$, there are
\begin{equation}\label{eqn:groups}
\frac{n!}{\displaystyle k_{211}!k_{111}!k_{011}!k_{200}!k_{100}!k_{000}!k_{210}!k_{110}!k_{010}!k_{201}!k_{101}!k_{001}!}
\end{equation}
different ways to partition $n$ vertices into twelve groups. Given a particular partition of the $n$ vertices, there are 
\begin{equation}\label{eqn:spoints}
{9 \choose 2}^{k_{211} + k_{210} + k_{201} + k_{200}}\cdot \left(9 \cdot 8 \cdot 7\right)^{k_{111} + k_{110} + k_{101} + k_{100}} \cdot \left( {9 \choose 2} {7 \choose 2} \right)^{k_{011} + k_{010} + k_{001} + k_{000}}
\end{equation}
different ways to assign special points. To illustrate the accuracy of this formula, consider a $k_{000}$ vertex. This vertex will consist of two out-points and seven in-points. Since two of the in-points in the first orientation are out-points in the other and because two of its out-points are in-points in the other orientation, the remaining five in-points that it consists of are in-points in both orientations. Similar reasoning is applied to the other eleven groups of vertices.

Next we will note that an (in,in)-point naturally pairs with an (out,out)-point. (Whenever we refer to say, (in,in)-points, the first ``in" refers to the first orientation and the second ``in" refers to the second orientation.) The number of (in,in)-points is
 \begin{align*}
&2k_{211} + 7k_{200} + k_{111} + k_{110} + k_{101} + 6 k_{100} + 5k_{000} + 2k_{010} + 2k_{001}\\
&= 2k_{211} + 2k_{200} + k_{111} - k_{110} - k_{101} + k_{100} + 2n + k - 2k_{210} - 2k_{201}.
\end{align*}
Likewise, the number of (out,out)-points is
\begin{equation}\label{eqn:pairs1}
2k_{211} + 2k_{200} + k_{111} - k_{110} - k_{101} + k_{100} + 2n + k - 2k_{210} - 2k_{201}.
\end{equation}
The remaining points will be (in,out)-points and (out,in)-points. The number of (in,out)-points is:
\begin{align}\label{eqn:pairs2}
\frac{\displaystyle 9n - 2(2k_{211} + 2k_{200} + k_{111} - k_{110} - k_{101} + k_{100} + 2n + k - 2k_{210} - 2k_{201})}{2} \nonumber \\
= 7n/2 - (2k_{211} + 2k_{200} + k_{111} - k_{110} - k_{101} + k_{100} + k - 2k_{210} - 2k_{201})
\end{align}
which also happens to be the number of (out,in)-points.

The set of (in,in)-points needs to be paired with the set of (out,out)-points. Similarly, (in,out)-points need to be paired with (out,in)-points. Hence, the product of (\ref{eqn:groups}), (\ref{eqn:spoints}), and the factorials of both (\ref{eqn:pairs1}) and (\ref{eqn:pairs2}) reflect the number of configurations. Dividing this product by the total number of matchings provides an estimate of the second-moment, which is the expected number of pairs of valid orientations. 

The second-moment is estimated as follows:
\begin{equation}\label{eq6}
\E[Y(Y-1)] = \sum\limits_{k \in I}r(\textbf{z})g(\textbf{z})\exp(nf(\textbf{z}))
\end{equation} 
where $\textbf{z} = \textbf{z}(\textbf{k}) = \textbf{k}/n$. The product within the sum is merely the result of applying Stirling's approximation $(s! \sim \sqrt{2 \pi n} (s/e)^s)$. The error term, the polynomial component, and the exponential component are represented by $r(\textbf{z})$, $g(\textbf{z})$, and $\exp(n f(\textbf{z}))$ respectively and defined below. Note that since $\sqrt{2 \pi n} ( \frac{n}{e})^n e^{\left( \frac{1}{12n} - \frac{1}{360n^3} \right)} < n! < \sqrt{2 \pi n} ( \frac{n}{e})^n e^{\left( \frac{1}{12n} \right)}$, $r(\textbf{z}) = O(1)$ for all $\textbf{z}$ and $r(\textbf{z}) \sim 1$ if $\textbf{z}$ is bounded away from the boundary of
\begin{align*}
J := \{ &(z,z_{211},z_{111}, z_{210}, z_{110},z_{201},z_{101},z_{200},z_{100}) \in  \mathbb{R}_{0}^{9} \bigm| z \leq 1/2, z_{211} + z_{111} \leq z, z_{200} \\
&+ z_{100} \leq z, z_{210} + z_{110} \leq 1/2 - z, z_{201} + z_{101} \leq 1/2 - z \}.
\end{align*}

By letting $b := 2z_{211} + 2z_{200} + z_{111} - z_{110} - z_{101} + z_{100} + 2 + z - 2 z_{210} - 2z_{201}$ and letting $h(x) := x \log x$ with the stipulation that $h(0) = 0$ we get 
\begin{align}\label{eqg}
g = &\left( b(9 - 2b)^{1/2} \right) / \bigm( \sqrt{128 (\pi n)^{9}} ( z_{211} z_{111} (z - z_{211} -z_{111}) z_{200} z_{100}(z - z_{200} - z_{100}) z_{210} z_{110} \nonumber \\
 &(1 - 2z - 2z_{210}  - 2 z_{110})  z_{201} z_{101} (1 - 2z - 2z_{201} -2z_{101}))^{1/2} \bigm).
\end{align}

The exponential contribution from (\ref{eqn:groups}) is
\begin{align*}
f_1 &= \log n - \log e\\
f_2 &= -h(z_{211}) -h(z_{111}) -h(z - z_{211} - z_{111}) - h(z_{200}) - h(z_{100}) - h(z - z_{200} - z_{100}) \\
& - h(z_{210})  - h(z_{110}) - h(1/2 -z - z_{210} - z_{110}) - h(z_{201}) - h(z_{101}) \\
& - h(1/2 - z - z_{201} - z_{101}) - \log n + \log e
\end{align*}
which correspond to the numerator and denominator, respectively. The exponential contribution from (\ref{eqn:spoints}) is
\begin{align*}
f_3 = &\log 9 + (z_{111} + z_{100} + z_{110} + z_{101}) \log 8 + (1 - z_{211} - z_{200} - z_{210} - z_{201}) \log 7 \\
& + (1 - z_{111} - z_{100} - z_{110} - z_{101})\log 4 \\
& + (1 - z_{211} - z_{200} -z_{210} - z_{201} - z_{111} - z_{100} - z_{110} - z_{101})\log 3.
\end{align*}
The exponential contribution from the factorials of (\ref{eqn:pairs1}) and  (\ref{eqn:pairs2}) is:
\begin{equation*}
f_4 =  h(9/2-b) + h(b) + 9/2(\log n -\log e).
\end{equation*}
The exponential contribution from dividing by the total number of matchings is
\begin{equation*}
f_5 = -h(9) + 9/2 \log2 + h(9/2) - 9/2\log n + 9/2 \log e.
\end{equation*}

Summing $f = f_1 + f_2+ f_3+ f_4 +f_5$ together we get
\begin{align}\label{eqf}
f =& h(b) + h(9/2 - b) +\log( 9 \cdot 4 \cdot 7 \cdot 3 ) - 9/2 \log9  - (z_{211} +z_{200} + z_{210} + z_{201} )\log7 \nonumber \\
& - (z_{211} + z_{200} + z_{210} + z_{200} + z_{111} + z_{100} + z_{110} + z_{101})\log3 \nonumber \\
& + (z_{111} + z_{100} + z_{110} + z_{101})\log 2 - h(z_{211}) - h(z_{111})\nonumber \\ 
&- h(z - z_{211} - z_{111}) - h(z_{200}) - h(z_{100}) - h(z - z_{200} - z_{100}) - h(z_{210}) - h(z_{110}) \nonumber \\
& - h(1/2 - z - z_{210} - z_{110}) - h(z_{201}) -h(z_{101}) - h(1/2 - z -z_{201} -z_{101}).
\end{align}
Since $h(0)=0$, $f$ is continuous towards the boundary of $J$, which allows us to explore the question of whether or not $f$ has a maximum. We calculate the partial derivatives of $f$ to identify stationary points in the interior of $J$. An algebraic manipulation tool such as Maple can be utilized towards this end.\footnote{The worksheets can be found on-line: \url{https://math.torontomu.ca/~pralat/}.}

We find that 
\begin{align*}
\frac{\delta f}{\delta z_{211}} &= \log \Bigm( \frac{\displaystyle b^2\cdot (z-z_{211} -z_{111})}{\displaystyle 7 \cdot 3 \cdot (9/2 - b)^2 \cdot z_{211}} \Bigm),\\
\frac{\delta f}{\delta z_{111}} &= \log \Bigm( \frac{\displaystyle  2 \cdot b \cdot  (z - z_{211} - z_{111})}{\displaystyle 3 \cdot  (9/2 - b) \cdot  z_{111}} \Bigm).
\end{align*}
Setting these two partial derivatives equal to one another and re-arranging gives us
$$
\frac{z_{111} }{z_{211}} = 14 \left( \frac{9/2 - b}{b}\right).
$$
Likewise, the partial derivatives for $z_{200}$ and $z_{100}$ are  
\begin{align*}
\frac{\delta f}{\delta z_{200}} &= \log \Bigm( \frac{\displaystyle b^2\cdot (z-z_{200} -z_{100})}{\displaystyle 7 \cdot 3 \cdot (9/2 - b)^2 \cdot z_{200}} \Bigm),\\
\frac{\delta f}{\delta z_{100}} &= \log \Bigm( \frac{\displaystyle  2 \cdot b \cdot  (z - z_{200} - z_{100})}{\displaystyle 3 \cdot  (9/2 - b) \cdot  z_{100}} \Bigm).
\end{align*}
Setting these two partial derivatives equal to one another and re-arranging gives us
$$
\frac{z_{100} }{z_{200}} = 14 \left( \frac{9/2 - b}{b}\right).
$$ 
Since $14 \left( \frac{9/2 - b}{b}\right) = \frac{z_{111} }{z_{211}} = \frac{z_{100} }{z_{200}}$, it follows that $z_{211}z_{100}=z_{200}z_{111}$.

The partial derivatives of $f$ with respect to $z_{210}$ and $z_{110}$ are
\begin{align*}
\frac{\delta f}{\delta z_{210}} &= \log \Bigm( \frac{\displaystyle (9/2 - b)^2\cdot (1/2 - z-z_{210} -z_{110})}{\displaystyle 7 \cdot 3 \cdot b^2 \cdot z_{210}} \Bigm),\\
\frac{\delta f}{\delta z_{110}} &= \log \Bigm( \frac{\displaystyle 2 \cdot (9/2 - b)\cdot (1/2 - z-z_{210} -z_{110})}{\displaystyle 7 \cdot 3 \cdot b \cdot z_{110}} \Bigm).
\end{align*}
The partial derivatives of $f$ with respect to $z_{201}$ and $z_{101}$ are
\begin{align*}
\frac{\delta f}{\delta z_{201}} &= \log \Bigm( \frac{\displaystyle (9/2 - b)^2\cdot (1/2 - z-z_{201} -z_{101})}{\displaystyle 7 \cdot 3 \cdot b^2 \cdot z_{201}} \Bigm),\\
\frac{\delta f}{\delta z_{101}} &= \log \Bigm( \frac{\displaystyle 2 \cdot (9/2 - b)\cdot (1/2 - z-z_{201} -z_{101})}{\displaystyle 7 \cdot 3 \cdot b \cdot z_{101}} \Bigm).
\end{align*}
Repeating the same procedure with regards to $\frac{\delta f}{\delta z_{210}}$,$ \frac{\delta f}{\delta z_{110}}$,$\frac{\delta f}{\delta z_{201}}$, and $\frac{\delta f}{\delta z_{101}}$ tells us that 
$$
z_{210}z_{101}=z_{201}z_{110}. 
$$

Setting $\frac{\delta f}{\delta z_{211}}=0$ gives $P_{211}=0$, where
$$
P_{211} =  b^2\cdot (z-z_{211} -z_{111}) - 7 \cdot 3 \cdot (9/2 - b)^2 \cdot z_{211}.
$$
If we define $P_{111}$, $P_{200}$, etc.\ similarly, and obtain $P$ from $\frac{\delta f}{\delta z}$, we obtain nine polynomials such that any local maximum in the interior of $I$ must be a common zero of all polynomials. We denote the resultant of two polynomials $X$ and $Y$ with respect to a variable $x$ as $\mathcal{R}(X,Y,x)$. When $X = Y = 0$ it is necessary that $ \mathcal{R}(X,Y,x) =0 $. Taking the resultant of $P_{211}$ and $P_{200}$ with respect to $z_{101}$ gives us
\begin{align*}
\mathcal{R}(P_{211},P_{200},z_{101}) &= 46294416(z z_{200} - z  z_{211} +z_{211}z_{100} -  z_{200}z_{111})^2\\
&= 46294416  ((z-z_{111})z_{200}-(z-z_{100})z_{211})^2.
\end{align*}
We can use the fact that $z_{211}z_{100}=z_{200}z_{111}$ and $\mathcal{R}(P_{211},P_{200},z_{101})=0$ to demonstrate that $z_{211}=z_{200}$, which then implies that $z_{111}=z_{100}$.

Calculating the resultant of $P_{201}$ and $P_{210}$ with respect to $z_{100}$ gives us
\begin{align*}
\mathcal{R}(P_{210},P_{201},z_{100}) &= 185177664(2z z_{201} - 2z z_{210} - 2z_{210}z_{101} + 2 z_{201} z_{110} - z_{201} + z_{210}),\\
&= 740710656((-1/2 + z +z_{110})z_{201} - (-1/2 + z +z_{101})z_{210})^2.
\end{align*}
We can use the fact that $z_{210}z_{101}=z_{201}z_{110}$ and $\mathcal{R}(P_{210},P_{201},z_{100})=0$ to demonstrate that $z_{210}=z_{201}$, which then implies that $z_{110}=z_{101}$.

We will now consider the nature of possible critical points. Let $z_{111} = z_{100}=c_1$, $z_{211} = z_{200}=c_2$, $z_{110} = z_{101}=c_3$, and $z_{210} = z_{201}=c_4$. 
Setting  $\frac{\delta f}{\delta z_{211}} = \frac{\delta f}{\delta z_{200}} = 0$ implies that 
$$
\left( \frac{9/2 - b}{b} \right)^2 = \frac{z - c_1 - c_2}{21 c_2}.
$$ After squaring $\frac{\delta f}{ \delta z_{100} } = \frac{\delta f}{ \delta z_{111} } = 0$, it follows that 
$$
\left( \frac{9/2 - b}{b}\right)^2 = \left( \frac{ 2(z - c_1 - c_2 ) }{3 c_1 }\right)^2.$$
Since 
$$\left( \frac{9/2 - b}{b} \right)^2 = \frac{z - c_1 - c_2}{21 c_2} = \left( \frac{ 2(z - c_1 - c_2 ) }{3 c_1 }\right)^2,$$
we have $z = \frac{3 c_1^2}{28 c_2} +c_1 + c_2$.

From setting $\frac{\delta f}{\delta z_{210}} = \frac{\delta f}{\delta z_{110}}$ we get $$
\frac{9/2 - b}{b} =  \frac{14 c_4}{c_3}.
$$ After setting $\frac{\delta f}{\delta z_{110}} = \frac{\delta f}{\delta z_{101} } = 0$ we get 
$$\frac{9/2 - b}{b} = \frac{3c_3}{2 (1/2 - z - c_4 -c_3)}.$$ 
Together these imply that 
$$
\left(\frac{9/2 - b}{b}\right)^2 = \frac{21 c_4}{(1/2 - z - c_3 - c_4)}.
$$ After squaring $\frac{\delta f}{\delta z_{110}} = \frac{\delta f}{\delta z_{101}}$, it follows that 
$$\left( \frac{9/2 - b}{b}\right)^2 = \left( \frac{3 c_3}{2(1/2 - z - c_3 - c_4)} \right)^2.$$ Since 
$$
\left( \frac{9/2 - b}{b}\right)^2 = \left( \frac{3 c_3}{2(1/2 - z - c_3 - c_4)} \right)^2 = \frac{21 c_4}{(1/2 - z - c_3 - c_4)},$$ it follows that $1/2 - z = \frac{3 c_3^2}{28 c_4} + c_3 + c_4$.

We now consider another parametrization of our variables that will hold when we are at a critical point. This parametrization is invaluable for finding the specific values of $b$ and $\textbf{z}$ that form a critical point for $f$. When $\frac{\delta f}{ \delta z_{111}} = \frac{\delta f}{ \delta z_{100}} = 0$, it follows that 
$$c_1 = 14 \left( \frac{9/2 - b}{b} \right)c_2.$$
 Likewise when $\frac{\delta f}{ \delta z_{211}} = \frac{\delta f}{ \delta z_{200}} = 0$, it follows that 
 $$b^2(z-c_1-c_2) = 21 \left(\frac{9}{2}-b \right)^2c_2, $$  which implies that 
$$z = \left( 21 \left(  \frac{9/2 - b}{b} \right)^2 + 14 \left( \frac{9/2 - b}{b} \right) + 1 \right) c_2.$$

When $\frac{\delta f}{\delta z_{110}} = \frac{\delta f}{\delta z_{101}} = 0$, 
$$c_3 = 14 \left( \frac{b}{9/2 - b} \right) c_4.
$$ Setting $\frac{\delta f }{\delta z_{210}} = \frac{\delta f }{\delta z_{201}} = 0$ implies that 
$$\left( 9/2 - b \right)^2 \left( \frac{1}{2} - z - c_3 - c_4 \right) = 21b^2c_4$$ which in turn implies that 
$$\frac{1}{2} - z = \left(  21 \left( \frac{b}{9/2 - b}\right)^2 + 14\left(\frac{b}{9/2 - b} \right) + 1 \right)c_4.$$ Since 
$$z = \left( 21 \left( \frac{9/2 - b}{b} \right)^2 + 14 \left( \frac{9/2 - b}{b} \right) + 1 \right) c_2,$$ 
it follows that
$$c_4 = \frac{1/2 - \left(21 \left( \frac{9/2 - b}{b} \right)^2 + 14 \left( \frac{9/2 - b}{b} \right) + 1 \right) c_2}{  21 \left( \frac{b}{9/2 - b}\right)^2 + 14\left(\frac{b}{9/2 - b} \right) + 1 }.$$ $c_3$ can be stated as 
$$c_3 = 14 \left(\frac{b}{9/2 - b}\right) \frac{1/2 - \left( 21 \left( \frac{9/2 - b}{b} \right)^2 + 14 \left( \frac{9/2 - b}{b} \right) + 1 \right) c_2}{  21 \left( \frac{b}{9/2 - b}\right)^2 + 14\left(\frac{b}{9/2 - b} \right) + 1 }.$$

After substituting these values for $c_1$, $c_2$, $c_3$, and $c_4$ into our partials, we are left with the task of discovering when our partials equal one another. It was found that for $\frac{\delta f }{\delta z_{211}} = \frac{\delta f }{\delta z_{200}} = \frac{\delta f }{\delta z_{111}} = \frac{\delta f }{\delta z_{100}}$, 
$$c_2^{*}(b) = - \frac{4 b^3 (32 b^2 +104 b - 171)}{5120 b^4 - 46080 b^3 + 285120 b^2 - 816480 b - 688905}$$ when $0 < b < 9/2$.

After substituting our values for $c_1$, $c_2$, $c_3$, and $c_4$ into $\frac{\delta f }{\delta z}$, we attain an expression for $\frac{\delta f }{\delta z}$ as a function of $b$ and $c_2$. We can substitute $c_2^{*}(b)$ for $c_2$ into our re-parameterized $\frac{\delta f }{\delta z}$ and then check which, if any, values of $b$ satisfy $\frac{\delta f }{\delta z}=0$. These values of $b$ are as follows: $b = b_1 \approx 0.8065779289$, $b =b_2 = 9/4$, and $b = b_3 \approx 3.693422071$. 

It now remains to see whether these values of $b$ satisfy 
$$0 < c_2 = c_2^{*}(b)< \frac{1/2}{ \left(21 \left( \frac{9/2 - b}{b} \right)^2 + 14 \left( \frac{9/2 - b}{b} \right) + 1 \right)}.$$ 
 Since $c_2^{*}(b_1) \approx -0.0001175309606 < 0$, $c_2^{*}(b_1)$ is not within the feasible region of $J$. 

The associated upper-bound 
$$c < \frac{1/2}{ \left(21 \left( \frac{9/2 - b}{b} \right)^2 + 14 \left( \frac{9/2 - b}{b} \right) + 1 \right)}$$ 
for $b_3$ is approximately 0.09883651395. Since $c_2^{*}(b_3) \approx 0.1105793451$, $c_2^{*}(b_3)$ is also not within the feasible region $J$.

We are now left with considering $b=b_2$, which corresponds to $c_2 = 1/144$. The upper-bound for this value is $1/72$, meaning that it is within the feasible region of $f$. This value and the corresponding values for $c_1, c_3, c_4$, and $z$ imply that $$ \hat{\textbf{z}}  = (1/4,1/144,7/72, 1/144,7/72,1/144,7/72,1/144,7/72)$$ is the sole critical point of $f$ in the interior of $J$. Note that 
$$
f(\hat{\textbf{z}} ) = \log \frac{756 \cdot \frac{9}{4}^{9/2} \cdot 2^{7/18} \cdot 12^{1/18}  }{3^{113/12} \cdot 7^{1/36} \cdot \frac{7}{72} ^{7/18} \cdot \frac{7}{48}^{7/12}} \approx 2.315007612.
$$

\subsection{Analyzing the Boundary Points}\label{subsec:boundary}

We will now consider candidate points to be global maxima which reside on the boundary of $J$. We will first consider any point on the boundary at which $0 < z < 1/2$. In this case, $0 \leq z_{211} + z_{111} \leq z$, $0 \leq z_{200} + z_{100} \leq z$, $0 \leq z_{110} + z_{210} \leq 1/2 - z$, and $0 \leq z_{101} + z_{201} \leq 1/2 - z$. 

When $z_{111} \rightarrow 0$, $\frac{\delta f}{\delta z_{111}}$ is dominated by the $z_{111}$ term on its denominator, which tends to $-\infty$. Likewise, when $z_{211} \rightarrow 0$, $\frac{\delta f}{\delta z_{211}}$ is dominated by the $z_{211}$ term on its denominator, which also tends to $-\infty$.

When $z_{211} + z_{111} \rightarrow z$, the $z - z_{211} - z_{111}$ term in the numerator of $\frac{\delta f}{\delta z_{111}}$ will tend to $\infty$. Likewise, the $z - z_{211} - z_{111}$ term in the numerator of $\frac{\delta f}{\delta z_{211}}$ will also tend to $\infty$. Thus, there is no global maximum at these particular points.

The same arguments also apply for $z_{200}$ and $z_{100}$, ensuring that no global maximum is reached at the points where either or both of them approach 0 and where their sum approaches $z$.

When $z_{110} \rightarrow 0$, $\frac{\delta f}{\delta z_{110}}$ is dominated by the $z_{110}$ term in its denominator, which tends to $-\infty$. Likewise, when $z_{210} \rightarrow 0$, $\frac{\delta f}{\delta z_{210}}$ is dominated by a $z_{210}$ term in its denominator, which also tends to $-\infty$. 

When $z_{210} + z_{110} \rightarrow 1/2 - z$, the $1/2 - z - z_{210} - z_{110}$ term in the numerator of $\frac{\delta f}{\delta z_{110}}$ will tend to $\infty$. Likewise, the $z - z_{210} - z_{110}$ term in the numerator of $\frac{\delta f}{\delta z_{210}}$ will also tend to $\infty$. Thus, there is no global maximum at these particular points.

The same arguments also apply for $z_{201}$ and $z_{101}$, ensuring that no global maximum is reached at the points where either or both of them approach 0 and where their sum approaches $1/2 - z$. 

It now remains to consider the points where $z = 0$ and $z = 1/2$. We first consider $z = 0$.

\medskip

\noindent \textbf{Case 1: $z = 0$}

When $z = 0$, $z_{211} = z_{111} = z_{200} = z_{100} = 0$ as well. We can thus state $f$ solely in terms of the remaining variables. Substituting $z = z_{211} = z_{111} = z_{200} = z_{100} = 0$ gives the function $\bar{f}(z_{210},z_{110},z_{201},z_{101})$ with domain $[0,1/2]^4$. 

Setting $\frac{\delta \bar{f}} { \delta z_{210}} = \frac{\delta \bar{f}}{ \delta z_{110}}$ gives us
\begin{equation*}
\frac{z_{210}}{z_{110}} = -\frac{(5 + 2 z_{110} + 2 z_{101} + 4 z_{210} + 4z_{201})}{28(2 - z_{110} - z_{101} - 2z_{210} - 2z_{201})}.
\end{equation*}
Likewise, setting $\frac{ \delta \bar{f} }{ \delta z_{201} }= \frac{\delta \bar{f}} { \delta z_{101}}$ gives us
$$
\frac{z_{201}}{z_{101}} = -\frac{(5 + 2 z_{110} + 2 z_{101} + 4 z_{210} + 4z_{201})}{28(2 - z_{110} - z_{101} - 2z_{210} - 2z_{201})}.
$$
Together these imply that $z_{210} z_{101} = z_{201} z_{110}$. 

Setting $\frac{\delta \bar{f}}{\delta z_{110}} = \frac{\delta \bar{f}}{\delta z_{101}} $ implies

$$
\frac{z_{110}}{1 - 2z_{210} - 2z_{110}} = \frac{z_{101}}{1 - 2z_{201} - 2z_{101}} 
$$
which implies that $z_{110} = z_{101}$.

Setting $\frac{\delta \bar{f}}{\delta z_{210}} = \frac{\delta \bar{f}}{\delta z_{201}} $ implies
$$
\frac{z_{210}}{1 - 2z_{210} - 2z_{110}} = \frac{z_{201}}{1 - 2z_{201} - 2z_{101}} 
$$
which implies that $z_{210} = z_{201}$. We are thus well-justified in considering $\bar{f}(z_{110},z_{210},z_{110},z_{210})$ in our search for critical points on this particular part of the boundary.

We now consider the segment $0 \leq z_{110} + z_{210} \leq 1/2$. The derivative of $\bar{f}(z_{110},z_{210},z_{110},z_{210}) $ with respect to $z_{110}$ is
$$
\frac{\delta \bar{f}(z_{110},z_{210},z_{110},z_{210}) } {\delta z_{110}} = \frac{(5 + 4z_{110} + 8 z_{210})(-1 + 2z_{210} + 2 z_{110})}{12(-1 + z_{110} + 2z_{210})z_{110}}
$$
which is not $0$ on the interior of the segment. It is only $0$ when $z_{110}+z_{210}=1/2$. 
$$
\frac{\delta \bar{f}(z_{110},z_{210},z_{110},z_{210}) } {\delta z_{210}} = -\frac{(5 + 4z_{110} + 8 z_{210})^2(-1 + 2z_{210} + 2 z_{110})}{168(-1 + z_{110} + 2z_{210})^2z_{210}}
$$
is also only $0$ when $z_{110}+z_{210}=1/2$. Differentiating $\bar{f}(z_{110},1/2-z_{110},z_{110},1/2-z_{110})$ with respect to $z_{110}$ produces
\begin{align*}
\frac{\bar{f}(z_{110},1/2-z_{110},z_{110},1/2-z_{110})}{\delta z_{110}}=&2 \log(2z_{110})-2\log(9 - 4z_{110}) +4 \log 2 + 2 \log 7 \\
& -2\log z_{110} +2 \log(1/2 - z_{110})
\end{align*}
which is $0$ when $z_{110} = 19/52$. Since 
\begin{align*}
& f(0,0,0,1/2-19/52,19/52,1/2-19/52,19/52,0,0) &=& \log \frac{\frac{133}{26} ^{19/26} \cdot \frac{98}{13}^{49/13}}{3^7 \cdot \frac{19}{52}^{19/26} \cdot 2^{27/26} \frac{7}{52}^{ 7/26} } \\
& \qquad \approx 1.672261141 &<& f(\hat{\textbf{z}} ),
\end{align*}
this point cannot be a global maximum.

\medskip

\noindent \textbf{Case 2: $z = 1/2$}

What follows for $z = 0$ is very similar to what was argued for $z= 0$. When $z=1/2$, $z_{210} = z_{110} = z_{201} = z_{101} = 0$ as well. We can thus state $f$ solely in terms of the remaining variables. Substituting $z_{210} = z_{110} = z_{201} = z_{101} = 0$ and $z = 1/2$ gives the function $\bar{f}(z_{211},z_{111},z_{200},z_{100})$ with domain $[0,1/2]^4$.

Setting $\frac{\delta \bar{f}} { \delta z_{211}} = \frac{\delta \bar{f}}{ \delta z_{111}}$ gives us
\begin{equation*}
\frac{z_{211}}{z_{111}} = -\frac{(5 + 2 z_{111} + 2 z_{100} + 4 z_{200} + 4z_{211})}{28(2 - z_{111} - z_{100} - 2z_{211} - 2z_{200})}.
\end{equation*}
Likewise, setting $\frac{\delta \bar{f}} { \delta z_{200}} = \frac{\delta \bar{f}}{ \delta z_{100}}$ gives us
\begin{equation*}
\frac{z_{200}}{z_{100}} = -\frac{(5 + 2 z_{111} + 2 z_{100} + 4 z_{200} + 4z_{211})}{28(2 - z_{111} - z_{100} - 2z_{211} - 2z_{200})}.
\end{equation*}
Together these imply that $z_{211}z_{100}=z_{200}z_{111}$.
Setting $\frac{\delta \bar{f}}{\delta z_{111}} = \frac{\delta \bar{f}}{\delta z_{100}} $ implies that
$$
\frac{z_{111}}{1 - 2z_{211} - 2z_{111}} = \frac{z_{100}}{1 - 2z_{200} - 2z_{100}},
$$
which in turn implies that $z_{111}=z_{100}$.
Setting $\frac{\delta \bar{f}}{\delta z_{211}} = \frac{\delta \bar{f}}{\delta z_{200}} $ implies that
$$
\frac{z_{211}}{1 - 2z_{211} - 2z_{111}} = \frac{z_{200}}{1 - 2z_{200} - 2z_{100}},
$$
which in turn implies that $z_{211}=z_{200}$. 

We are thus well-justified in considering $\bar{f}(z_{111},z_{211},z_{111},z_{211})$ in our search for critical points on this particular part of the boundary.
We now consider the segment $0 \leq z_{111} + z_{211} \leq 1/2$. The derivative of $\bar{f}(z_{111},z_{211},z_{111},z_{211}) $ with respect to $z_{111}$ is
$$
\frac{\delta \bar{f}(z_{111},z_{211},z_{111},z_{211}) } {\delta z_{111}} = \frac{(5 + 4z_{111} + 8 z_{211})(-1 + 2z_{211} + 2 z_{111})}{12(-1 + z_{111} + 2z_{211})z_{111}}
$$
which is only zero when $z_{211} = z_{111} = 1/2$.

Differentiating $\bar{f}(z_{111},1/2-z_{111},z_{111},1/2-z_{111})$ with respect to $z_{111}$ produces
\begin{align*}
\frac{\bar{f}(z_{111},1/2-z_{111},z_{111},1/2-z_{111})}{\delta z_{111}}=&2 \log(2z_{111})-2\log(9 - 4z_{111}) +4 \log 2 + 2 \log 7 \\
& -2\log z_{111} +2 \log(1/2 - z_{111})
\end{align*}
which is $0$ when $z_{111} = 19/52$. Since 
\begin{align*}
& f(1/2,1/2-19/52,19/52,0,0,0,0,1/2-19/52,19/52) &=& \log \frac{\frac{133}{26} ^{19/26} \cdot \frac{98}{13}^{49/13}}{3^7 \cdot \frac{19}{52}^{19/26} \cdot 2^{27/26} \frac{7}{52}^{ 7/26} } \\
& &\approx& 1.672261141 < f(\hat{\textbf{z}} ) , 
\end{align*}
this point cannot be a maximum.

The Hessian matrix of $f(\hat{\textbf{z}})$ is
$$B = \frac{1}{63}
\left[
\begin{array}{ccccccccc}
-1672 & 544 & 488 & 544 & 488 & -544 & -488 & -544 & -488 \\
544 & -9280 & -320 & 224 & 112 & -224 & -112 & -224 & -112 \\
488 & -320 & -1024 & 112 & 56 & -112 & -56 & -112 & -56\\
544 & 224 & 112 & -9280 & -320 & -224 & -112 & -224 & -112\\
488 & 112 & 56 & -320 & -1024 & -112 & -56 & -112 & -56\\
-544 & -224 & -112 & -224 & -112 & -9280 & -320 & 224 & 112\\
-488 & -112 & -56 & -112 & -56 & -320 & -1024 & 112 & 56\\
-544 & -224 & -112 & -224 & -112 & 224 & 112 & -9280 & -320\\
-488 & -112 & -56 & -112 & -56 & 112 & 56 & -320 & -1024\\
\end{array}
\right].
$$
To calculate the signs of the eigenvalues of $B$, we can calculate the eigenvalues of $B^*:=63 B$.
The characteristic polynomial of $B^*$ is 
$$(x^2 + 10584x + 10077696)^3  (x^3 + 11136x^2 + 21055680x + 3072577536).$$
The eigenvalues of $B^*$ are 
{\small
\begin{align*}
\lambda_1 &= \lambda_2 = \lambda_3 = -108\sqrt{1537} - 5292 \approx -9526.09588932514,\\
\lambda_4 &= \lambda_5 = \lambda_6 = 108\sqrt{1537} - 5292\approx -1057.90411067487,\\
\lambda_7 &=  -8\sqrt{3} \sqrt{105631}\sin((\pi - \arctan((9\sqrt{5834559781407})/26571068))/3) \\
&- 8\sqrt{105631}\cos((\pi - \arctan((9\sqrt{5834559781407})/26571068))/3) - 3712 \approx -8776.89694570535,\\
\lambda_8 &=  8\sqrt{3}\sqrt{105631}\sin((\pi - \arctan((9\sqrt{5834559781407})/26571068))/3) \\
& - 8\sqrt{105631} \cos((\pi - \arctan((9\sqrt{5834559781407})/26571068))/3) - 3712 \approx -2199.97604519778,\\
\lambda_9 &=  16\sqrt{105631}\cos((\pi - \arctan((9\sqrt{5834559781407})/26571068))/3) - 3712\approx -159.127009096879.
\end{align*}} Since all of these values are negative the Hessian matrix at $\hat{\textbf{z}}$ is negative definite, implying that $\hat{\textbf{z}}$ is a local maximum. The determinant of $B$ is found to be $\frac{-23665185138564661248}{117649}$.

We will now approximate $\E [Y^2]$ using a method reminiscent of Lemma 6.3 in \cite{greenski}, to which the reader is referred to for a more detailed perspective on what follows. To the end of approximating $\E [Y^2]$, we will set $z = 1/4 + y$, $z_{2jk} = 1/144 + y_{2jk}$, and $z_{1jk} = 7/72 + y_{1jk}$. Setting $\textbf{y} = (y,y_{211},y_{111},y_{200},y_{100},y_{210},y_{110},y_{201},y_{101})$ and using Taylor's Theorem to expand around $ \hat{\textbf{z}}$ gives us 
$$
f = \log(81/8) + \frac{\textbf{y}B \textbf{y}^T}{2} + O(x^3).
$$
Let $x= x(\textbf{z}) = ||\textbf{y}||$ with $||\cdot||$ denoting the $L^2$ norm. Let $J_0 = \{\textbf{z} : x \leq n^{-2/5} \}$. For $\textbf{z} \in J_0$ we have $r(\textbf{z})g(\textbf{z}) \sim g(\hat{\textbf{z}})$ and $x^3 = O(n^{-6/5})$ since $J_0$ is bounded away from the boundary of $J$ and $\textbf{z} \in J_0$ are close to $\hat{\textbf{z}}$.

\par Taking all of this together gives us
\begin{align}
\sum\limits_{\textbf{k}:\textbf{k}/n \in J_0}r(\textbf{z})g(\textbf{z})\exp(n f(\textbf{z})) &= \sum\limits_{\textbf{k}:\textbf{k}/n \in J_0}r(\textbf{z})g(\textbf{z})\exp(n f(\hat{\textbf{z}} + n \textbf{y} B \textbf{y}^T + O(x^{-3})))\nonumber, \\
&\sim \left( \frac{81}{8} \right)^n g(\hat{\textbf{z}}) e^{O(n^{-1/5})} \sum\limits_{\textbf{k}:\textbf{k}/n \in J_0} e^{n \textbf{y} \frac{B}{2} \textbf{y}^T	}\nonumber ,\\
&\sim \left( \frac{81}{8} \right)^n g(\hat{\textbf{z}}) \sum\limits_{\textbf{k}:\textbf{k}/n \in J_0} e^{n \textbf{y} \frac{B}{2} \textbf{y}^T	} \label{eq7}.
\end{align}
Dividing $\sum\limits_{\textbf{k}:\textbf{k}/n \in J_0} e^{n \textbf{y} \frac{B}{2} \textbf{y}^T	}$ by $n^9$ gives us a Riemann sum as $n \to \infty$. By changing $\textbf{k}$ to $\textbf{z}=\textbf{k}/n$ and defining $\textbf{w} = \sqrt{n} \textbf{y}$ to rescale the region of integration we get
$$
\frac{1}{n^9} \sum\limits_{\textbf{k}:\textbf{k}/n \in J_0} e^{n \textbf{y} \frac{B}{2} \textbf{y}^T	}= \frac{1}{n^9} \sum\limits_{\textbf{k}:\textbf{k}/n \in J_0} e^{\textbf{w} \frac{B}{2} \textbf{w}^T	} = \int_{J_0} e^{\textbf{w} \frac{B}{2} \textbf{w}^T}d^9\textbf{w}. 
$$
Since the side-length of the scaled region is $\sqrt{n} \cdot n^{-2/5} = n^{1/10}$ which goes to $\infty$ as $n \to \infty$ we have
$$
\int_{J_0} e^{\textbf{w} \frac{B}{2} \textbf{w}^T}d^9\textbf{w} \sim \int_{\mathbb{R}^9} e^{\textbf{w} \frac{B}{2} \textbf{w}^T}d^9\textbf{w} = \int_{\mathbb{R}^9} e^{n \textbf{y} \frac{B}{2} \textbf{y}^T}d^9\textbf{y}.
$$

After diagonalizing and using the Gaussian integral $\int_{-\infty}^\infty e^{-x^2} = \sqrt{\pi}$ we get 
$$
\int_{\mathbb{R}^9} e^{n \textbf{y} \frac{B}{2} \textbf{y}^T}d^9\textbf{y} = \frac{1}{n^9} \cdot \sqrt{ \frac{n^{9}\cdot \pi^{9}}{|\det \frac{B}{2}|}} = \frac{1}{n^9} \cdot \sqrt{ \frac{n^{9} \cdot 2^9 \cdot \pi^{9}}{|\det B|}} = \frac{2^{9/2} \pi^{9/2} }{n^{9/2} \sqrt{|\det B|}}
$$
which leaves us with
\begin{equation}\label{domeq}
\left( \frac{81}{8} \right)^n g(\hat{\textbf{z}}) \sum\limits_{\textbf{k}:\textbf{k}/n \in J_0} e^{n \textbf{y} \frac{B}{2} \textbf{y}^T	} = \left( \frac{81}{8} \right)^n g(\hat{\textbf{z}}) \frac{(2 \pi n)^{9/2}}{ \sqrt{|\det B|}}.
\end{equation}
It now remains to consider points in $J\backslash J_0$.

On the boundary of $J_0$ the value of $f$ is $f(\hat{\textbf{z}} ) - \Omega(n^{-4/5}) $. Recall that $B$ is negative definite (meaning $\textbf{y} B \textbf{y}^T <0$ for all $\textbf{y} \in \mathbb{R}^9$). Furthermore, $f$ is a fixed function that is independent of $n$ and has a single global maximum at $\hat{\textbf{z}}$. This implies that
$$
\max \limits_{\textbf{z} \in J \backslash J_0} f(\textbf{z}) = f(\hat{\textbf{z}} )- \Omega(n^{-4/5}) = \log ( 81/8) - \Omega(n^{-4/5}).
$$

Since $r$ and $g$ are polynomially bounded it follows that the terms in the summation (\ref{eq6}) for which $\textbf{k}/n \in J/J_0$ are bounded by $(81/8)^n e^{-\Omega(n^{1/5})}$. Since there are polynomially many such terms their sum is $(81/8)^n e^{-\Omega(n^{1/5})}$.

It follows, then, that
\begin{align}
\E[Y(Y-1)] &= \sum\limits_{\textbf{k}:\textbf{k}/n \in J}r(\textbf{z})g(\textbf{z})\exp(n f(\textbf{z}))\nonumber ,\\ 
 &= \sum\limits_{\textbf{k}:\textbf{k}/n \in J_0}r(\textbf{z})g(\textbf{z})\exp(n f(\textbf{z})) + \sum\limits_{\textbf{k}:\textbf{k}/n \in J \backslash J_0}r(\textbf{z})g(\textbf{z})\exp(n f(\textbf{z}))\nonumber,\\
&\sim \left( \frac{81}{8} \right)^n g(\hat{\textbf{z}}) \frac{2^{9/2} \pi^{9/2} }{n^{9/2} \sqrt{|\det B|}} + (81/8)^n e^{-\Omega(n^{1/5})} \nonumber,\\
&\sim \left( \frac{81}{8} \right)^n g(\hat{\textbf{z}}) \frac{2^{9/2} \pi^{9/2} }{n^{9/2} \sqrt{|\det B|}}\nonumber,\\
&\sim \left( \frac{81}{8} \right)^n \frac{81}{7}.
\end{align}
Thus we have verified Lemma~\ref{lemma:second} and Corollary~\ref{cor:second} follows:
$$
\frac{\E[Y(Y-1)]}{\E[Y]^2} = \frac{\left( \frac{81}{8} \right)^n \frac{81}{7}}{ 9\left(\frac{81}{8} \right)^{n}}=  \frac{9}{7}.
$$
\section{Joint Factorial Moments}\label{sec:joint}
To complete our proof of Theorem \ref{thm:jaegp} we introduce random variables $X_1, X_2, \ldots, X_k$ for $k \geq 1$ which represent the number of cycles of length $k$ in $\mathcal{P}(n,9)$. Here condition (a) of the Small Subgraph Conditioning Method (Theorem~\ref{thm:sscm}) follows from a well known result by Bollob\'as (see Corollary 2.19 in \cite{bolo2}) stipulates that $X_1,X_2, \ldots, X_k$ are asymptotically independent Poisson distributions with mean
$$
\E [X_k] \longrightarrow \lambda_k := \frac{8^k}{2k}.
$$

To satisfy the condition (b) for the small subgraph conditioning method and thereby prove Theorem \ref{thm:jaegp}, we will need to show that there is a constant $\mu_k$ for each $k \geq 1$ such that 
$$
\frac{\E [Y X_k]}{\E Y} \longrightarrow \mu_k
$$
and that the joint factorial moments satisfy 
$$
\frac{\E [Y [X_1]_{j_1} [X_2]_{j_2} \ldots [X_k]_{j_k}]}{\E Y} \longrightarrow \prod\limits_{i=1}^k \mu_i^{j_i}
$$
for any fixed $j_{1}, j_2, \ldots, j_n$, where $[x]_k$ is the falling factorial $x! / (x-k)!$. To estimate $\E [Y X_k]$ we will count the mean number of triples $(P,C,O)$, where $P$ represents a pairing, $C$ represents a $k$-cycle, and $O$ represents an orientation by dividing the number of such pairings by the total number of matchings $|\mathcal{P}(9,n)| = M(9n)$.
	
There are $\frac{(9 \cdot 8)^k}{2k} \cdot \frac{n!}{(n-k)!}$ different ways to choose vertices that make up the cycle. There are $i$ vertices with in-degree 2, $b$ of which are in-vertices. Since the sum of in-degrees in the cycle must equal the number of out-degrees in the cycle, there are $i$ vertices in the cycle with out-degree 2, $c$ of which are out-vertices. This implies that $k-2i$ vertices in the cycle have an in-degree and out-degree of 1 in the context of the cycle. We let $a_i = 2 {k \choose 2i}$ refer to the number of orientations. The number of choices of in-vertices among $i$-vertices with in-degree 2 is ${i \choose b}$. Likewise, there are ${i \choose c}$ possible choices of out-vertices among $i$ vertices with out-degree 2. Turning our attention to the vertices outside of the cycle, we see that there are ${n-2i \choose n/2 - i -b +c}$ ways to select the in- and out-vertices that are not in $C$. 
	
We now consider the distribution of special points between our orientations. There are ${7 \choose 2}^{i-b}$ ways to select special points in vertices that have an in-degree of 2 in the cycle and are out-vertices in the graph. Likewise, there are ${7 \choose 2}^{i-c}$ ways to select special points in vertices that have an out-degree of 2 in the cycle and which are in-vertices in the graph. The special points for vertices with in-degree 2 in $C$ and which are in-vertices and for vertices with out-degree 2 in $C$ and which are out-vertices do not need to be identified. There are $7^{k-2i}$ ways to select special points among vertices in $C$ which have an in-degree of 1 and an out-degree of 1. There are ${9 \choose 2}^{n-k}$ ways to select special points for vertices outside of $C$. Given that $9n/2 -k$ must be in-points, there are $(9n/2 - k)!$ different ways to pair in-points with out-points. We are now ready to verify Lemma~\ref{lemma:cycles}. 

\begin{proof}[Proof of Lemma~\ref{lemma:cycles}.]
We calculate $\frac{\E [Y X_k ]}{\E[Y]}$ as follows: 	
\begin{eqnarray*}
\frac{\E [Y X_k ]}{\E[Y]} &\sim& \frac{\sum\limits_{i=0}^{\lfloor k/2 \rfloor} \frac{[n]_k (9 \cdot 8)^k}{2k} a_i {9 \choose 2}^{n-k} (9n/2 - k)! 7^{k-2i} \sum\limits_{b=0}^i \sum\limits_{c=0}^i {i \choose b}{i \choose c} {n-2i \choose n/2 - i - b + c} {7 \choose 2}^{i-b} {7 \choose 2}^{i-c}  }{ { n \choose n/2} {9 \choose 2}^n (9n/2)!},\\
&\sim &  \frac{\sum\limits_{i=0}^{\lfloor k/2 \rfloor} \frac{[n]_k (9 \cdot 8)^k}{2k} a_i 7^{k-2i} \sum\limits_{b=0}^i\sum\limits_{c=0}^i {i \choose b}{i \choose c} {n-2i \choose n/2 - i - b + c} {7 \choose 2}^{i-b}{7 \choose 2}^{i-c}}{{n \choose n/2} (9 \cdot 4)^k (9n/2)^k},\\
&\sim & \frac{ \sum \limits_{i=0}^{\lfloor k/2 \rfloor} \frac{n^k}{2k} a_i 2^k 7^{k-2i} \frac{(n-2i)!}{n!} \sum \limits_{b=0}^i {i \choose b} (7 \cdot 3)^b \sum\limits_{c=0}^i {i \choose c} (7 \cdot 3)^c \frac{n/2 !}{(n/2-i+b-c)!}\frac{n/2!}{(n/2-i-b+c)!}}{(\frac{9n}{2})^k}.\\
\end{eqnarray*}
Since $n \to \infty$ and $k$ is finite, Stirling's approximation implies that $\frac{(n-2i)!}{n!} \sim \left(\frac{n}{e} \right)^{-2i}$ and\\  $\frac{n/2 !}{(n/2-i+b-c)!}\frac{n/2!}{(n/2-i-b+c)!} \sim \left( \frac{n}{2e}\right)^{2i}$. We therefore get
\begin{eqnarray*}
\frac{\E [Y X_k ]}{\E[Y]} &\sim& \frac{ \sum \limits_{i=0}^{\lfloor k/2 \rfloor} \frac{n^k}{2k} a_i 2^k 7^{k-2i} \frac{(n-2i)!}{n!} \sum \limits_{b=0}^i {i \choose b} (7 \cdot 3)^b \sum\limits_{c=0}^i {i \choose c} (7 \cdot 3)^c \frac{n/2 !}{(n/2-i+b-c)!}\frac{n/2!}{(n/2-i-b+c)!}}{(\frac{9n}{2})^k},\\
&\sim& \sum \limits_{i=0}^{\lfloor k/2 \rfloor} \left( \frac{28}{9}\right)^k \frac{a_i}{2k} \sum\limits_{b=0}^i {i \choose b} 21^b  \sum \limits_{c=0}^i {i \choose c} 21^c \left(\frac{n}{e} \right)^{-2i} \left( \frac{n}{2e}\right)^{2i} .\\
\end{eqnarray*}

Since $\sum\limits_{b=0}^{i} {i \choose b}21^b = \sum\limits_{c=0}^{i} {i \choose c}21^c = (1+21)^i$ we get the following:
\begin{eqnarray*}
\frac{\E [Y X_k ]}{\E Y} &\sim & \sum \limits_{i=0}^{\lfloor k/2 \rfloor} \frac{1}{2k} \left( \frac{28}{9} \right)^k a_i \left(\frac{1}{7} \right)^{2i} 22^{2i} \left(\frac{1}{2} \right)^{2i},\\
&\sim& \frac{1}{2k}\left( \frac{28}{9} \right)^k \sum\limits_{i=0}^{\lfloor k/2 \rfloor} 2 {k \choose 2i} \left( \frac{11}{7}\right)^{2i}.
\end{eqnarray*}
This can be completed by defining $q(x):= 2(1+x)^k = \sum\limits_{i=0}^k 2 {k \choose i} \cdot i$. We can consider our particular sum to be the even terms of the summation $\sum\limits_{i=0}^{\lfloor k/2 \rfloor} 2 {k \choose 2i} \left( \frac{11}{7}\right)^{2i}$. Since the exponential term will be negative for odd $i$ we can simplify our expression as follows:
\begin{eqnarray*}
\sum\limits_{i=0}^{\lfloor k/2 \rfloor} 2 {k \choose 2i} \left( \frac{11}{7}\right)^{2i} &=&\frac{1}{2}\left( q\left(\frac{11}{7}\right) + q\left(-\frac{11}{7}\right)\right).\\
\end{eqnarray*}
Returning to our calculation of $\frac{\E [Y X_k]}{\E Y}$ gives us
\begin{eqnarray*}
\frac{\E [Y X_k ]}{\E Y} &=& \frac{1}{2k} \left( \frac{28}{9} \right)^k \left( \left(\frac{18}{7} \right)^k +  \left(\frac{-4}{7}\right)^k \right)
= \frac{1}{2k} \left( 8^k + \left( -\frac{16}{9} \right)^k \right)\\
&=& \frac{8^k}{2k} \left( 1 + \left( -\frac{2}{9} \right)^k \right)= \lambda_k \left( 1 + \left( -\frac{2}{9} \right)^k \right).\\
\end{eqnarray*}
We conclude that $\mu_k  = \lambda_k \left( 1 + \left( -\frac{2}{9} \right)^k \right)$.
\end{proof}

Since Corollary 2.19 from \cite{bolo2} stipulates that the cycles are asymptotically independent, Lemma~1 by Janson from~\cite{janson1995random} implies that for every finite sequence $j_1, j_2, \ldots, j_k$ of non-negative integers, the following holds:
$$
\frac{\E(Y [X_1]_{j_1} \cdots [X_k]_{j_k})}{\E Y} \sim \prod\limits_{i=1}^k \mu_i^{j_i} = \prod\limits_{i=1}^k \left(  \lambda_i(1+\delta_i ) \right)^{j_i}
$$
where $\delta_i = \frac{\mu_i}{\lambda_i} - 1$ to satisfy condition (c). For the $i$ that we concerned ourselves with above $\delta_i = -\left(\frac{2}{9}\right)^i$.

\bibliography{bibjaeg}

\end{document}